\documentclass{elsart}
\usepackage{amsfonts}
\usepackage{amssymb}
\usepackage{graphicx}

\begin{document}

\begin{verbatim}\end{verbatim}\vspace{2.5cm}

\begin{frontmatter}
\title{ Helly $\mathbf{EPT}$ graphs on bounded degree trees: forbidden induced subgraphs and efficient recognition}
\author[laplata]{L. Alc\'on}
\ead{liliana@mate.unlp.edu.ar}
\author[laplata,coni]{M. Gutierrez}
\ead{marisa@mate.unlp.edu.ar}
\author[laplata,coni]{M. P. Mazzoleni}
\ead{pia@mate.unlp.edu.ar}

\address[laplata]{Departamento de Matem\'atica,
Universidad Nacional de La Plata,  CC 172, (1900) La Plata,
Argentina}
\address[coni]{CONICET}

\begin{abstract}
The edge intersection graph of a family of paths in  host tree is
called an $EPT$ graph. When the host tree has maximum degree  $h$,
we say that $G$ belongs to the class $[h,2,2]$. If, in addition,
the family of paths satisfies the Helly property, then $G \in$
Helly $[h,2,2]$.
 The time complexity of the recognition of the classes
 $[h,2,2]$ inside the class $EPT$ is open for every $h> 4$.
In \cite{gol}, Golumbic et al. wonder if the only obstructions for
an $EPT$ graph belonging to $[h,2,2]$ are the chordless cycles
$C_n$ for $n> h$. In the present paper, we   give a negative
answer to that question,  we present a   family of
$EPT$ graphs 
which are forbidden induced subgraphs for the classes $[h,2,2]$.
Using them we obtain a total characterization by induced forbidden
subgraphs of the classes Helly $[h,2,2]$ for $h\geq 4$ inside the
class $EPT$. As a byproduct, we prove that Helly $EPT$$\cap
[h,2,2]=$ Helly $[h,2,2]$. Following the approach used in
\cite{Monma}, we characterize Helly $[h,2,2]$ graphs by their
atoms in the decomposition by clique separators. We give an
efficient algorithm to recognize Helly $[h,2,2]$
graphs.\end{abstract}\end{frontmatter}


\section{Introduction}

A graph $G$ is called $EPT$ (or $UE$) if it is the edge
intersection graph of a family of paths in a tree. $EPT$ graphs
are used in network applications, the problem of scheduling
undirected calls in a tree network is equivalent to the problem of
coloring an $EPT$ graph (see \cite{ver}).  The class of $EPT$
graphs was first investigated by Golumbic and Jamison
\cite{B3,B4}. In the last decades many papers were devoted to the
study of $EPT$ graphs and their generalizations, see
\cite{20years,GolTrenk, Spinrad}. In \cite{B6},  the class of
graphs that admit an $EPT$ representation on a host tree with
maximum degree $h$ is denoted by $[h,2,2]$. Clearly, $[2,2,2]$ is
the class of interval graphs. It is known that $[3,2,2]$ is
precisely the class of chordal $EPT$ graphs \cite{B6}, while
$[4,2,2]$ is the class of weakly chordal $EPT$ graphs \cite{gol4}.
Notice that the class of $EPT$ graphs is the union of the classes
$[h,2,2]$ for
 $h\geq 2$. A complete hierarchy of
related graph classes  emerging by imposing different restrictions
on the tree representation is published in \cite{gol}.

On the algorithmic side, the recognition and coloring problems
restricted to $EPT$ graphs are NP-complete, whereas the maximum
clique and maximum stable set problems are polynomially solvable.
See \cite{B3}.

The time complexity of the recognition of the classes
 $[h,2,2]$ inside the class $EPT$ is open for $h> 4$, and it is known
 to be polynomial time solvable for $h\in\{2,3,4\}$.
In \cite{gol} and \cite{gol4}, Golumbic et al. wonder if the only
obstructions for an $EPT$ graph belonging to $[h,2,2]$ are the
chordless cycles of size greater than $h$. In \cite{pia}, we give
a negative answer to this question and present a family of
forbidden induced subgraphs called prisms.

In this paper, we generalize the class of prisms and present a
wider family of $EPT$ graphs called $k$-gates which are forbidden
induced subgraphs for the classes $[h,2,2]$ when $h<k$.

A graph is Helly $EPT$ (or $UEH$) if it admits an $EPT$
representation using a path family that satisfies the Helly
property. In \cite{Monma}, Monma and Wei characterize $EPT$ and
Helly $EPT$ via decomposing the graph by clique separators and
prove that the latter class can be recognized efficiently. Finding
a characterization by forbidden induced subgraphs of $EPT$ and of
Helly $EPT$ graphs are long standing open problems.

 Helly $[h,2,2]$ is the class of  graphs that admit a Helly $EPT$
representation on a host tree with maximum degree $h$. Clearly,
 Helly $EPT$$\cap [h,2,2]\subseteq $ Helly $[h,2,2]$ but the equality not necessary holds.

We obtain a total characterization by induced forbidden subgraphs
of the class Helly $[h,2,2]$ inside the class $EPT$ using gates.
gates. As a byproduct, we prove that Helly $EPT$$\cap [h,2,2]=$
Helly $[h,2,2]$ which means that, in the way of way of getting a
Helly representation, it is not necessary to increase the maximum
degree of the host tree.

\begin{figure} \centering{
\includegraphics[width=10cm]{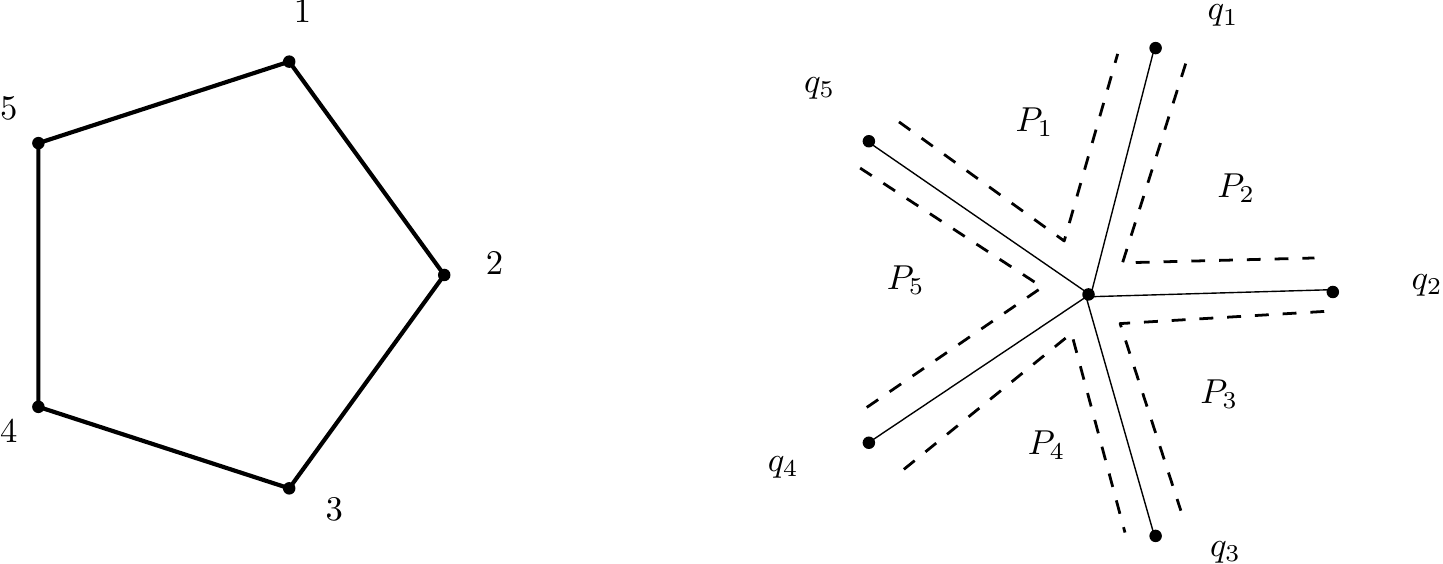}}
\caption{The cycle $C_5$ and its $EPT$ representation: a pie of
size $5$.}\label{f:ciclo}
\end{figure}

In addition, we characterize Helly $[h,2,2]$ graphs by their atoms
in the decomposition by clique separators. We give an efficient
algorithm to recognize Helly $[h,2,2]$ graphs.

 The paper is
organized as follows: in Section \ref{s:prelim}, we provide basic
definitions and known results. In Section \ref{d:gates}, we depict
the graphs named $k$-gates and focus on their  main properties; we
show that $k$-gates are Helly $EPT$ but do not admit an $EPT$
representation on a host tree with maximum degree less than $k$.
In Section \ref{s:HEPT_charac}, we show that  a Helly $EPT$ graph
$G$ belongs to the class Helly $[h,2,2]$ if and only if $G$ does
not have a $k$-gate as induced subgraph for any $k>h$. Finally, in
Section \ref{s:complex}, we  use the Monma and Wei decomposition
by clique separator to obtain an efficient algorithm for the
recognition of Helly $[h,2,2]$ graphs.

\section{Preliminaries and known results} \label{s:prelim}
%
In this paper all graphs are finite and simple. Given a graph $G$,
$V(G)$ and $E(G)$ denote the vertex set and the edge set of $G$,
respectively. An $\mathbf{EPT}$ \textbf{representation} of $G$ is
a pair $\langle \mathcal{P},T \rangle$ where $\mathcal{P}$ is a
family $(P_v)_{v\in V(G)}$ of subpaths of the \textbf{host tree}
$T$ satisfying that two vertices $v$ and $w$ of $G$ are adjacent
if and only if $E(P_v)\cap E(P_w) \neq \emptyset$. When the
maximum degree of  the host tree $T$ is $h$, the $EPT$
representation of $G$ is called an
$\mathbf{(h,2,2)}$\textbf{-representation} of $G$. The class of
graphs that admit an $(h,2,2)$-representation is denoted by
$\mathbf{[h,2,2]}$.

A \textbf{star} is any complete bipartite graph $K_{1,n}$. The
only vertex with degree grater than one is called  the
\textbf{center of the star}. The edges of a star are called
\textbf{spokes}. The star $K_{1,3}$ is named the \textbf{claw
graph}. We will say that a path $P:(v_1,...,v_l)$ contains a
vertex $v$ if $v=v_i$ for some $1\leq i \leq l$;  and that it
contains an edge $e$ if $e=v_iv_{i+1}$ for some $1\leq i \leq
l-1$.

Golumbic et al. introduced the notion of pie in order to describe
$EPT$ representations of chordless cycles. A \textbf{pie of size }
$\mathbf{k}$ in an $EPT$ representation $\langle \mathcal{P},T
\rangle$ is a star subgraph of $T$ with central vertex $q$ and
neighbors $q_1$,...,$q_{k}$ and a subfamily of  paths
$P_1,...,P_k$ of $ \mathcal{P}$ such that
$\{q_i,q,q_{i+1}\}\subseteq V(P_i)$, for $1 \leq i \leq k$
(addition is assumed to be module $n$). See Figure \ref{f:ciclo}.
\begin{thm} \cite{B3} \label{t:pie} Let
$\langle \mathcal{P},T \rangle$ be an $EPT$ representation of a
graph $G$. If $G$ contains a chordless cycle $C_k$ with $k\geq 4$,
then $\langle \mathcal{P},T \rangle$ contains a pie of size $k$
whose paths are in one-to-one correspondence with the vertices of
$C_k$.\end{thm}
A set family $(S_i)_{i\in I}$ satisfies the \textbf{ Helly
property} if any pairwise intersecting subfamily  $(S_i)_{i\in
I'}$ with $\emptyset \neq I'\subseteq I$ has non-empty total
intersection, i.e. $\bigcap_{i\in I'}S_i\neq \emptyset$. A graph
$G$ is \textbf{Helly} $\mathbf{EPT}$ if it admits an $EPT$
representation $\langle \mathcal{P},T \rangle$ such that the set
family $(E(P))_{P\in \mathcal{P}}$ satisfies the Helly property.
In an analogous way, we say that $G$ is \textbf{Helly}
$\mathbf{[h,2,2]}$  if it admits an $(h,2,2)$-representation
$\langle \mathcal{P},T \rangle$ such that the family $(E(P))_{P\in
\mathcal{P}}$ satisfies the Helly property. Clearly, Helly
$[h,2,2] \subseteq$ Helly $EPT \cap [h,2,2]$.

 A \textbf{complete set} of a graph $G$ is a subset of  $V(G)$
whose elements are pairwise adjacent. A \textbf{clique} is a
maximal (with respect to the inclusion relation) complete set.

\begin{figure} \label{f:cliques}\centering{
\includegraphics[width=10cm]{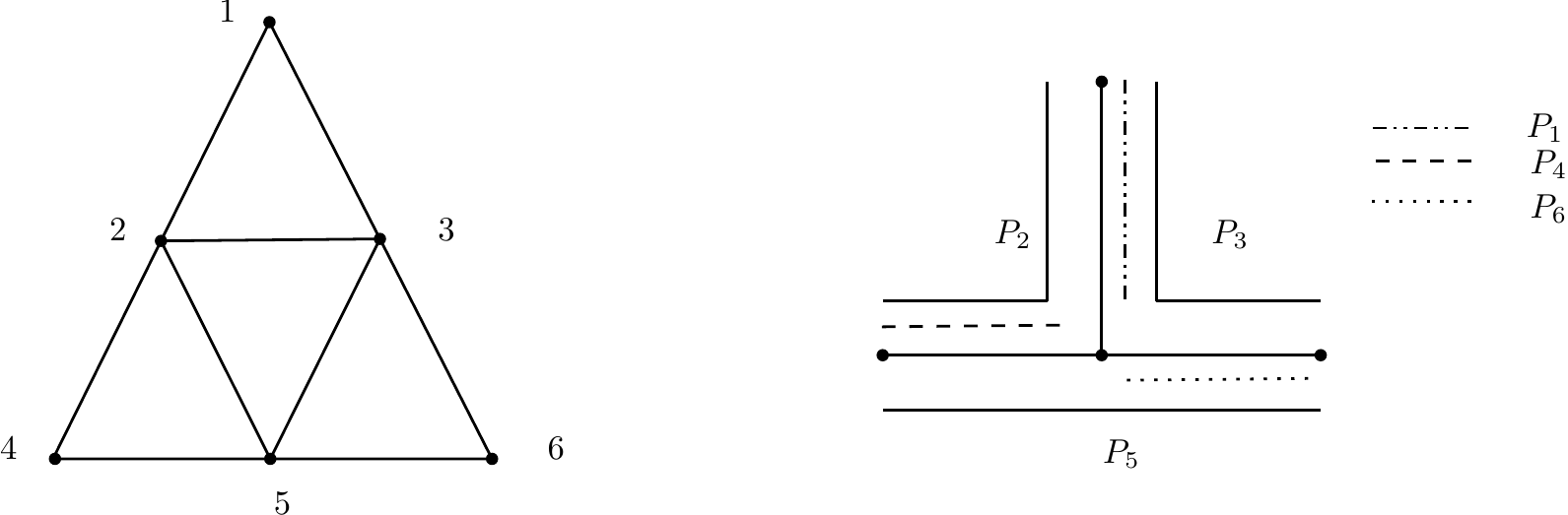}}
\caption{An $EPT$ representation of the sun $S_3$. In this
representation, the central triangle $\{2, 3, 5\}$ is a
claw-clique; the other three triangles are edge-cliques.}
\end{figure}

Given an $EPT$ representation $\langle (P_v)_{v\in V(G)},T
\rangle$  of $G$, for every edge $e$  of $T$, let $K_e$ be the
complete set $ \{v\in V(G): e \in E(P_v)\}$. For every claw $Y$ in
$T$, let $K_Y$  be the complete set $\{v\in V(G):  P_v$ contains
two spokes of $Y\}$.
%
\begin{thm}\cite{B4}
Let $\langle \mathcal{P},T \rangle$ be an $EPT$ representation of
$G$. If $C$ is a clique of $G$ then either there is an edge $e\in
E(T)$ such that $C=K_e$ or there is a claw $Y$ in $T$  such that
$C=K_Y$.
\end{thm}
%
In the former case, when there exists $e$ such that $C=K_e$, the
clique $C$ is called an \textbf{edge-clique}, otherwise $C$ is
called a \textbf{claw-clique}. See Figure \ref{f:cliques}. Notice
that the condition of being an edge-clique or a claw-clique
depends on the given representation. Clearly, in a Helly $EPT$
representation every clique is an edge-clique.
%
%
We say that three paths   of a given $EPT$ representation $\langle
\mathcal{P},T \rangle$ \textbf{form a claw} if there exists a claw
$Y$ of $T$ such that every pair of spokes of $Y$ is contained by
some of the paths. Clearly, there is claw-clique if and only if
three paths form a claw.

If $S \subseteq V(G)$ then $G-S$ denotes the graph induced in $G$
by $V(G)\setminus S$. When $S$ contains a unique vertex $v$, we
 write simply $G-v$.
\section{Gates and multipies} \label{s:gates}
%
A clear corollary of Theorem \ref{t:pie} is that every chordless
cycle $C_k$ with $k>  h\geq 3$ is an obstruction for the class
$[h,2,2]$. In \cite{gol}, Golumbic et al. wonder if besides cycles
there are other $EPT$ forbidden induced subgraphs for this class.
In \cite{pia}, answering negatively the previous question, we
described for every $h> 4$ an $EPT$ graph $F_h$ which has no
induced cycles of size $k$ for every $k>h$, but it does not admit
an $EPT$ representation on a host tree with maximum degree less
than or equal to $h$. The graphs introduced in the following
definition generalize the graphs $F_h$. In Section
\ref{s:HEPT_charac}, we obtain a total characterization of Helly
$[h,2,2]$ graphs using them.

\begin{figure} \label{f:tranquera} \centering{}
\includegraphics[width=12cm]{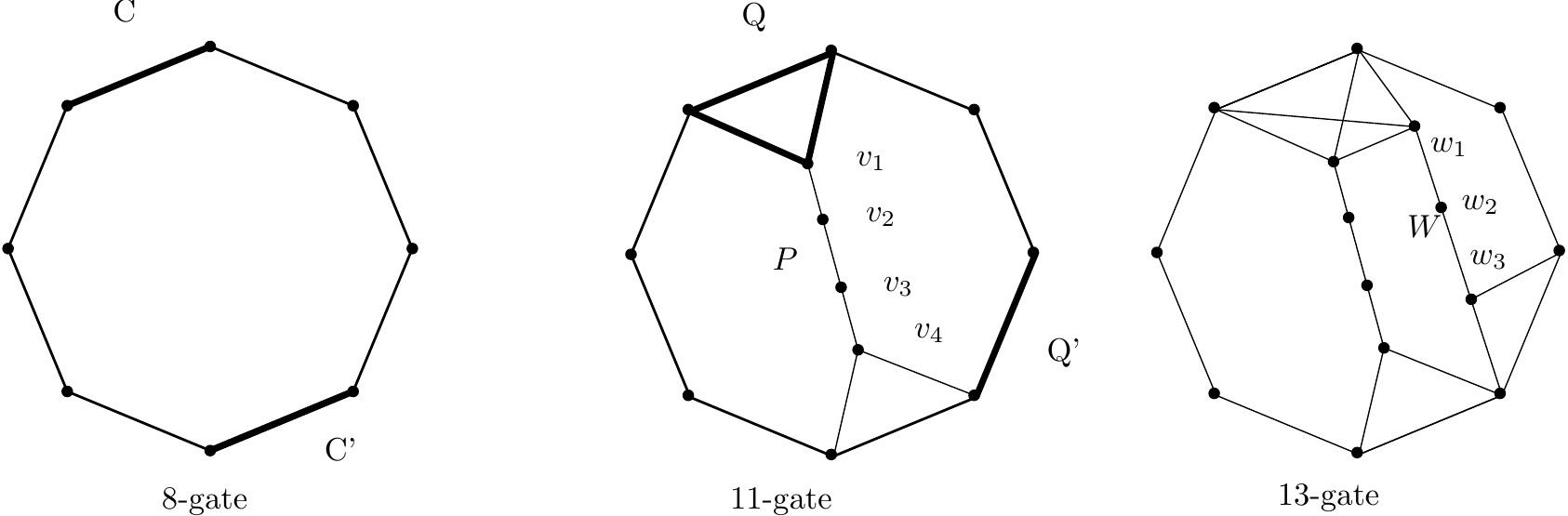}
\caption{Some examples of gates. From left to right, the second
gate is obtained from the first using the bold cliques $C$ and
$C'$ and the path $P:(v_1,v_2,v_3,v_4)$. The  third gate is
obtained from the second using the bold cliques $Q$ and $Q'$ and
the path $W:(w_1,w_2,w_3)$. }
\end{figure}

\begin{figure}
\centering{
  \includegraphics[width=10cm]{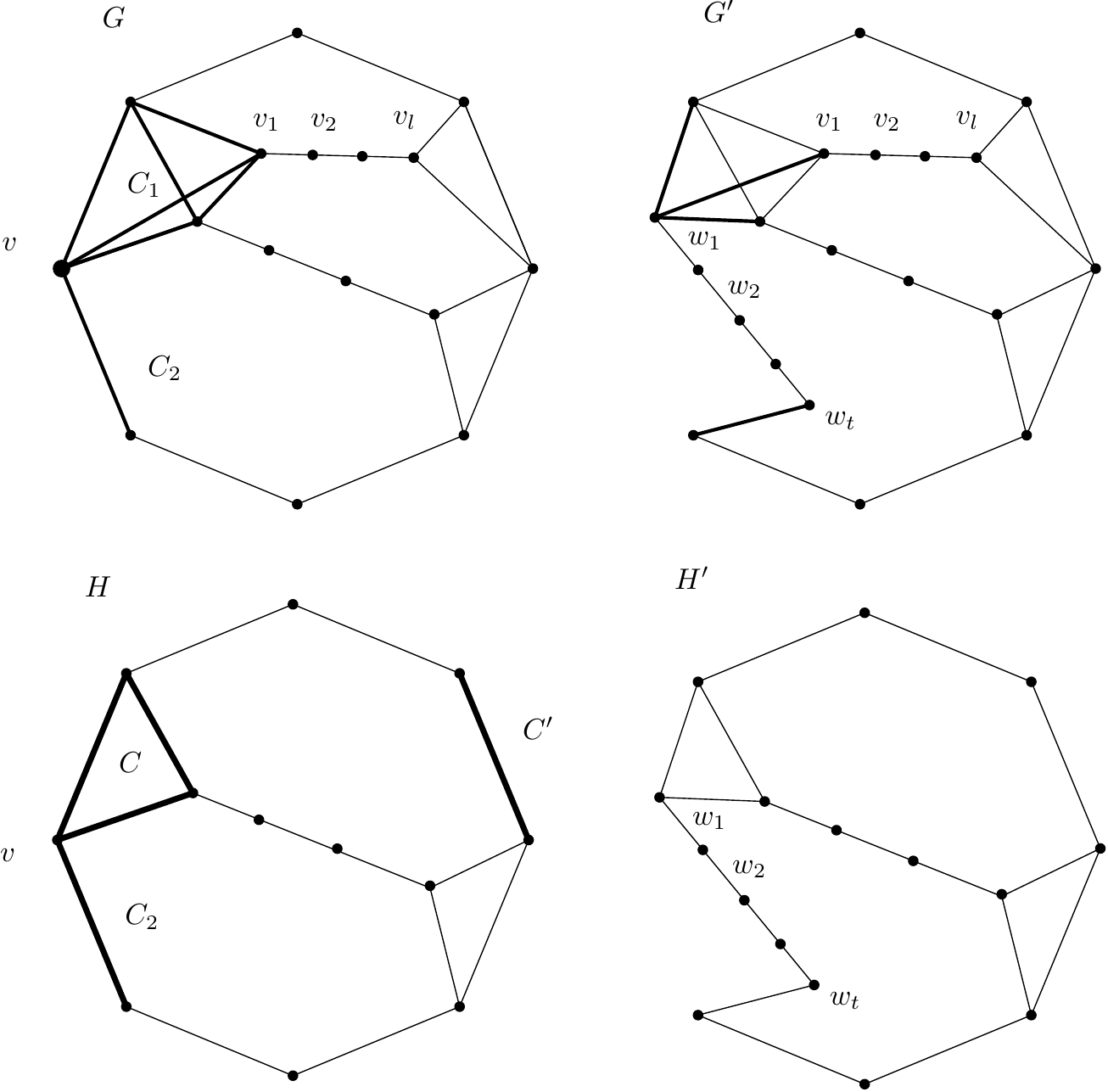}}
  \caption{An example following the proof of Lemma \ref{l:nodisj}.}\label{f:nodisj}
\end{figure}

We say that two graphs $G$ and $G'$ are \textbf{disjoint} if
$V(G)\cap V(G')=\emptyset$. The \textbf{union} of $G$ and $G'$ is
the graph $H$ with $V(H)=V(G) \cup V(G')$ and $E(H)= E(G) \cup
E(G')$.
%
%
\begin{defn}\label{d:gates} The following graphs are called \textbf{gates}.
\begin{itemize}
\item Every chordlees cycle $C_n$ with  $n\geq 4$ is a gate;\
\item If $G$ is a gate, $C$ and $C'$ are disjoint cliques of $G$,
and $P:( v_1,..,v_l)$ with
 $l\geq 2$ is a chordless path disjoint from $G$, then
the union of $G$ and $P$ plus all edges between $v_1$ and the
vertices of $C$, and all edges between $v_l$ and the vertices of
$C'$ is a gate;\
 \item There are no more gates.\end{itemize}
If the number of cliques of a gate $G$ is $k$ then we say that $G$
is a \textbf{k-gate}.\end{defn}
%
 In Figure \ref{f:tranquera} we offer some examples of gates.
\begin{lem}\label{l:gateEPT}
If $G$ is a $k$-gate then $G\in$ Helly $[k,2,2]$. Furthermore, $G$
admits a Helly $(k,2,2)$-representation  on a host tree that is a
star.
\end{lem}
\begin{pf} We proceed by induction.
   Clearly the statement holds for $C_k$.

If $G$ is not a cycle, then $G$ is obtained from an $m$-gate $H$
 using disjoint cliques $C$ and $C'$ of $H$ and a path
 $P:(v_1,v_2,..,v_{l})$ with $l\geq2$ disjoint from $H$.  Notice that $m+(l-1)=k$. Let $\langle \mathcal{P},T
 \rangle$ be a Helly $(m,2,2)$-representation of $H$ with $T$ a
 star. We can assume that $T$ has $m$ spokes. Let $e$ and $e'$
 be  spokes of $T$ such that $C=K_e$ and $C'=K_{e'}$. Denote by $T'$
the star that is obtained by adding  $l-1$ spokes
$e_1,...,e_{l-1}$ to $T$. Let $P_{v_1}$ be the subpath of $T'$
defined by the edges $e$ and $e_1$. For $2\leq i \leq l-1$ let
$P_{v_i}$ be the subpath of $T'$ defined by the edges $e_{i-1}$
and $e_{i}$; and let $P_{v_l}$ the one defined by the edges
$e_{l-1}$ and $e'$.

Thus  $\langle \mathcal{P'},T'
 \rangle$ is a Helly $(k,2,2)$-representation of $G$, where $\mathcal{P'}$ is the
 family
 $\mathcal{P}$ plus the paths $P_{v_i}$ for $1 \leq i \leq l$.
\hfill$\Box$\end{pf}
%
\begin{lem}\label{l:unico} If $G$ is a gate and $v\in V(G)$, then $v$ belongs to
exactly two cliques of $G$. In addition, if $C_1$ and $C_2$ are
those cliques then $C_1\cap C_2=\{v\}$.\end{lem}
\begin{pf} We proceed by induction.
   Clearly the statement holds for chordless cycles.

Let $G$ be a gate obtained from another gate $H$, using disjoint
cliques $C$ and $C'$ of $H$ and a  chordless path $P:(
v_1,..,v_l)$ with $l\geq 2$  disjoint from $H$. Notice that the
cliques of $G$ are:
\begin{itemize}
    \item[]  the cliques of $H$ other than $C$ and $C'$;
 \item[] the cliques of $P$, i.e. $\{v_i,v_{i+1}\}$ for $1\leq i\leq l-1$;
     \item[]$C\cup \{v_1\}$; and
      \item[] $C'\cup \{v_l\}$.
\end{itemize}
The proof follows easily from the fact that $H$ satisfies the
statement.\hfill$\Box$\end{pf}
\begin{lem}\label{l:nodisj} Let $v$ be a vertex of a gate $G$,
 $C_1$ and $C_2$ cliques of $G$ such that $
 C_1\cap C_2=\{v\}$,
  and
$W:( w_1,..,w_t)$  a chordless path disjoint from $G$ with $t\geq
2$. Then, the graph $G'$ union of $G-v$ and $W$ plus all edges
between $w_1$ and the vertices of $C_1-\{v\}$ and all edges
between $w_t$ and the vertices of $C_2-\{v\}$ is a gate.\end{lem}
\begin{pf} We proceed by induction. Clearly
 the statement holds for chordless cycles.

Assume $G$ is a gate obtained from another gate $H$, using
disjoint cliques $C$ and $C'$ of $H$ and a chordless path $P:(
v_1,..,v_l)$ with $l\geq 2$  disjoint from $H$.\
 If $v$ is one of
the vertices  of $P$ then the proof is direct and simple.

If $v$ is a vertex of $C$ (see Figure \ref{f:nodisj}),  we can
assume that $C_1=C\cup \{v_1\}$ and $C_2$ is a clique of $G$
different from $C'\cup \{v_l\}$, which means that in $H$ the
vertex $v$ is the intersection between the cliques $C$ and $C_2$ .
Thus, by the inductive hypothesis, the graph $H'$ obtained from
the  union of $H-v$ and $W$ plus all edges between $w_1$ and the
vertices of $C-\{v\}=C_1-\{v,v_1\}$ and all edges between $w_t$
and the vertices of $C_2-\{v\}$ is a gate. Since the path $P$ is
disjoint from $H'$, and $(C_1-\{v_1,v\}) \cup\{w_1\}$ and $C'$ are
disjoint cliques of $H'$, thus, by the recursive definition of
gate, the union of $H'$ and $P$ plus all edges between $v_1$ and
the vertices of $(C_1-\{v_1,v\}) \cup\{w_1\}$, and all edges
between $v_l$ and the vertices of $C'$ is a gate. The proof
follows from the fact that this is the same graphs $G'$ depicted
in the statement of the theorem.

If $v$ is a vertex of $C'$ or if $v\in V(H)-(C\cup C')$ the proof
is analogous.\hfill$\Box$\end{pf}
%
 Golumbic and Jamison proved  that (see Theorem \ref{t:pie})
 chordless
cycles admit a \textit{unique} $EPT$ representation called pie. In
what follow, generalizing that result, we introduce the definition
of multipie and prove that  also gates  admit a \textit{unique}
$EPT$ representation.
%
\begin{defn} \label{d:multipie} A \textbf{multipie of size} $\mathbf{k}$ in an
$EPT$ representation $\langle \mathcal{P},T \rangle$ is a star
subgraph of $T$ with central vertex $q$ and neighbors
$q_1$,..,$q_k$ and a  subfamily $\mathcal{P'}$ of $\mathcal{P}$
such that:
 \begin{enumerate}
 \item if $P\in \mathcal{P'}$ then $|V(P)\cap \{q_1,q_2,..,q_k\}|=2$  (every path contains
 two spokes of the star);\
 \item if $i\neq j$ then $|\{P\in \mathcal{P'}:
\{q_i,q_j\}\subseteq V(P)\}|\leq 1$ (no two paths contain the same
two spokes);
 \item if $1\leq i\leq k$ then $|\{P\in \mathcal{P'}: \{q,q_i\}\subseteq V(P)\}|\geq 2$   (every spoke of the
 star is contained by at least two paths);\
 \item no three paths of $\mathcal{P'}$ form a claw.
\end{enumerate}\end{defn}
%
%
Observe  that every pie is a multipie. The following theorem
generalize Theorem \ref{t:pie}.
\begin{thm} \label{t:multipie} Let $\langle \mathcal{P},T \rangle$ be an $EPT$ representation
of $G$. If $G$ contains a
$k$-gate  then $\langle \mathcal{P},T \rangle$ contains a multipie
of size $k$   whose paths are in one-to-one correspondence with
the vertices of the gate.\end{thm}
%
\begin{pf} Let $\langle \mathcal{P},T \rangle$ be an $EPT$ representation of $G$ whit
$\mathcal{P}=(P_v)_{v\in V(G)}$.
 We can assume, without loss of generality,
 that $G$ is a $k$-gate. We proceed by induction.
If $G$ is a chordless cycle $C_k$ then, by Theorem \ref{t:pie},
$\langle \mathcal{P},T \rangle$ contains a pie of size $k$ and the
proof follows.

If $G$ is not a cycle, then $G$ is obtained from an $m$-gate $H$
 using disjoint cliques $C$ and $C'$ of $H$ and a path
 $P:(v_1,v_2,..,v_{l})$ with $l\geq 2$ disjoint from $H$.
 Notice that $m+(l-1)=k$.  By inductive
hypothesis, $\langle \mathcal{P},T \rangle$ contains a multipie of
size $m$ formed by  a star subgraph $S$ of $T$ and the  path
subfamily  $\mathcal{P'}=(P_v)_{v\in V(H)}$.

Let $S$ be the star with center $q$ and leaves $q_1,...,q_m$. By
condition $(4)$ in  Definition \ref{d:multipie}, no three paths of
$\mathcal{P'}$ form a claw, then there exists a spoke of $S$, say
$e_1=qq_1$, such that $C\subseteq K_{e_1}$; and there exists
another spoke, without loss of generality say $e_2=qq_2$, such
that $C'\subseteq K_{e_2}$. Even more, by condition $(2)$, $e_1$
and $e_2$ are the only spokes of $S$ satisfying the described
property.

Let $d$ be the minimum distance in $H$ between a vertex of $C$ and
a vertex of $C'$. Clearly, $d\geq 1$. Chose vertices $u\in C$ and
$u'\in C'$ such that the distance between them in $H$ is $d$. Let
$(u,u_1,u_2,...,u_{d-1},u')$ be a shortest path in $H$ between
 $u$ and $u'$. Notice that $u,u_1,u_2,...,u_{d-1},u',
v_l,v_{l-1},...v_2,v_1$ induce a cycle in $ G$ of size $d+l+1\geq
4$. By Theorem \ref{t:pie}, in $\langle \mathcal{P},T \rangle$
there is a pie corresponding to this cycle. Let $S'$ be the star
subgraph of $T$ used by this pie. Notice that the center of $S'$
must be the same vertex $q$ of $T$. Even more, since the vertex
$v_1$ of $P$ is adjacent to all vertices in $C$, the vertex $v_l$
is adjacent to all vertices in $C'$, and there are no other
adjacencies  between vertices of $P$ and $H$, we have that $S'$
has $l-1$ spokes that are not spokes of $S$. The remaining
$(d+l+1)-(l-1)=d+2$ spokes of $S'$ are also spokes of $S$.
Therefore the union of $S$ and $S'$ is a star subgraph of $T$ with
center $q$ and $m+l-1=k$ spokes. Now it is not difficult to check
that $\mathcal{P}$ forms a multipie around the star $S\cup S'$,
and the proof follows. \hfill $\Box$\end{pf}
%
\section{Forbidden induced subgraphs for Helly EPT graphs on bounded degree trees}
\label{s:HEPT_charac}
%

The goal of this section is Theorem \ref{t:bueno} below. We prove
that  gates are the only  subgraphs which force the use a host a
tree with large enough degree in every Helly $EPT$ representation
of a graph.
%
\begin{thm} \label{t:bueno} Let $G$ be a Helly $EPT$ graph and $h\geq 3$.
 Then, $G\notin$ Helly $[h,2,2]$
 if and only if there exists $k>h$
 such that $G$ has a $k$-gate as induced subgraph. \end{thm}
%
\begin{pf} We will prove the direct implication,
the converse  follows from Theorem \ref{t:multipie} and the fact
that Helly $[h,2,2]\subseteq [h,2,2]$.\

 Assume that $G$ is a Helly $EPT$ graph which does not admit
 a Helly $(h,2,2)$-representation. Let $d$ be the smallest positive integer such that  $G\in$ Helly $[d,2,2]$.
 Clearly, $d>h$. Let $\langle \mathcal{P},T \rangle$ be a Helly
 $(d,2,2)$-representation of $G$
 minimizing the number of vertices of the host tree  $T$ with degree $d$.
%
 \begin{claim}\label{palitos} We can assume that if $q\in V(T)$ is the end vertex of
     a path $P\in \mathcal{P}$ then $d_T(q)\leq 2$.\end{claim}
%
\begin{pf} If it is  not the case, by subdividing every edge of $T$ (and consequently
every edge of every path of $\mathcal{P}$)   in three
    parts, and after that shortening  every path of $\mathcal{P}$ by removing its two  end vertices,
     we obtain the desired representation.\hfill $\Box$\end{pf}

    Let $q_0\in V(T)$ be a vertex with degree $d$ and call
    $q_1,... ,q_d$ to its neighbors. Denote by $H$  the subgraph
    of $G$ induced by the vertices $v$ such that $q_0\in V(P_v)$.
%
%
    \begin{claim}\label{haytranquera}
      The subgraph $H$ contains an induced cycle of length at least 4.\end{claim}
%
\begin{pf}  Let $P=(v_1,...,v_l)$ be the longest induced path in $H$ and
assume, without loss of generality, that
$\{q_i,q_0,q_{i+1}\}\subseteq V(P_{v_i})$, for all $i:1,..,l$.
Notice that $2\leq l\leq d-1$.
%
%

    Suppose, in order to derive a contradiction, that
     every path of $\mathcal{P}$ containing $q_0q_{l+1}$
    also contains $q_0q_l$.     Then, we can modify (as explained below) the
    representation $\langle \mathcal{P},T\rangle$ to obtain a new
    Helly $(d,2,2)$-representation of $G$ on a host tree with less vertices of
    degree $d$, contrary to our assumption.
    Indeed, to obtain the new representation do:

    subdivide the edge $q_0q_l$ adding a new vertex $\widetilde{q}_l$;

    remove the edge $q_0q_{l+1}$ and do $q_{l+1}$ adjacent to
    $\widetilde{q}_l$;

    in the paths of $\mathcal{P}$ containing the edge
    $q_0q_{l+1}$, replace the vertex $q_0$ and the edges  $q_0q_{l+1}$ and $q_0q_l$ by the vertex
    $\widetilde{q}_l$ and the edges
    $\widetilde{q}_lq_{l+1}$ and $\widetilde{q}_lq_l$, respectively;

no other path is modified except for the fact of subdividing the
edge $q_0q_l$ if necessary.

     Therefore, there must exist $1\leq j\leq d$, $j\neq l, l+1$,  and a vertex $x$ of $H$   such that $\{q_j, q_0 ,q_{l+1}\}\subseteq
    V(P_x)$. Clearly, $x\notin V( P)$.

    If $j> l+1$, then $V(P)\cup \{x\}$
     induces a path of $H$ longer than $P$, which contradicts the election of $P$.

   If $j=l-1$, then $P_x$, $P_{v_{l-1}}$ and $P_{v_{l}}$ violate the Helly property,
    which contradicts the fact that this a Helly $EPT$
    representation of $G$.

 Thus $j\leq l-2$. This implies that  $H$ contains
 the cycle induced by the vertices $\{v_j,,...,v_{l-1},v_{l},x\}$,
  as we wanted to prove.   \hfill $\Box$\end{pf}

It follows from the previous Claim \ref{haytranquera}, that $H$
has at least an induced gate.  Let $R$ be a biggest induced gate
in $H$, say that $R$ is a $k$-gate, and assume without loss of
generality, that the multipie corresponding to the vertices of $R$
use the star with edges $\{q_0q_1,...,q_0q_k\}$ (see Lemma
\ref{t:multipie}).

 We will prove that $k=d$.  Since $d>h$,  the proof
follows.

Clearly, $k\leq d$. Suppose, in order to derive a contradiction,
that $k<d$.

 Since $G$ is
    connected there must exists a vertex $y$ such that the path $P_y$  uses one of the
    edges $q_0q_1$,...,$q_0q_k$ and an edge $q_0q_i$ for some $k< i \leq d$.
     Without loss of generality, we can assume
     that $\{q_k,q_0,q_{k+1}\} \subseteq V(P_y)$.

If all  paths containing the edge $q_0q_{k+1}$ also contain the
edge $q_0q_k$,  then (as we did before) we can modify  the
    representation $\langle \mathcal{P},T\rangle$ to obtain a new
    representation of $G$ on a host tree with fewer vertices of
    degree $d$, contrary to assumption.

    So, there exists a vertex $z$ and $j\neq k,k+1$ such that
    $\{q_j,q_0,q_{k+1}\} \subseteq V(P_z)$.
     Notice that $y$ and $z$ are adjacent and
    do not belong to the gate $R$.

    Assume, in order to derive a contradiction, that $j\leq k-1$. Let $C_k$
    and $C_j$ be the  cliques of $R$ corresponding to the edges $q_0q_k$
    and $q_0q_j$ of $T$, respectively. Notice that $C_k$ and $C_j$ are disjoint, otherwise
    $P_y$, $P_{z}$ and $P_{v}$  violate the Helly property, where $v$ is a vertex in the intersection.
    Using  cliques $C_k$, $C_j$ and the path $P:(y,z)$
    disjoint from $R$, we obtain a $(k+1)$-gate induced in $H$,
    which contradicts the election of $R$. Therefore, $j>k-1$. Since  $j\neq k, k+1$, say $j=k+2$.

    Denote by $A$  the set of vertices  $v\in V(H)$ such that
    $P_v$ contains an edge $q_0q_i$
    for some $i\leq k$ and an edge $q_0q_{i'}$
    for some $i'> k$. Notices that $y\in A$ and $z\not\in A$. Let $G_{z}$ be the
    connected component of $G-A$ containing the vertex $z$.

Clearly, if $v\in V(G_z)\cap V(H)$ then there exist $i$ and $i'$,
$k+1\leq i< i'\leq d$ such that $\{q_i,q_0,q_{i'}\}\subseteq
    V(P_v)$, thus,
    without loss of generality, we can assume that there exists $s$, with
    $k+2\leq s \leq d$, such that

    $$V(G_{z})\cap V(H)=\bigcup_{k+1\leq i < i' \leq s}\{v\in V(G):
     \{q_i,q_0,q_{i'}\}\subseteq V(P_v)\};$$
and for every  $k+1\leq i  \leq s$
    \begin{equation}\label{eq:1}\hbox{ there exists
     $v\in V(G_{z})\cap V(H)$ such that $q_0q_i\in
     E(P_v)$.}\end{equation}
%
 \begin{claim} \label{c:q0qk} If    $y'\in A$ and $P_{y'}$ contains an edge $q_0q_t$ with
    $ k+1 \leq t  \leq s$ then  $P_{y'}$ also contains the edge $q_0q_k$.
    \end{claim}
%
 \begin{pf} Assume, in order to derive a contradiction, that $q_0q_{j}\in E(P_{y'})$
with $1\leq j<k$. Since $q_0q_t\in E(P_{y'})$ and  $ k+1 \leq t
\leq s$, by (\ref{eq:1}),  there exists
 $z'\in V(G_{z})\cap V(H)$ adjacent to $y'$. We chose $z'$
 minimizing its distance to $z$ in $G_z$ (it could be $z'=z$).
 Let $P:(z,z_1,...,z')$ be a
 shortest $zz'$-path in $G_z$. It is clear that $y'$ es adjacent to no vertex
 of $P$ except $z'$. Notice also that
 $V(P)\cap V(R)=\emptyset$, and no vertex of $P$ is adjacent to a
 vertex of $R$. So we will deal with the following cases:
 \begin{itemize}
    \item[(a)] $y$ is adjacent to $y'$ (in this case $t$ must be equal to $k+1$).
    \item[(b)] $y$ is non adjacent to $y'$ but it is adjacent to some vertex of
    $P$ besides  $z$. Thus, it must be $z_1$ and
    $(y,z_1,...,z',y')$ is a chordless path.
    \item[(c)] $y$ is neither adjacent to $y'$ nor to a vertex of $P$
    besides $z$. Thus $(y,z,z_1,...,z',y')$ is a chordless path.
 \end{itemize}

 Let $C_j$ and $C_k$ be the cliques of $R$
 corresponding to de edges $q_0q_j$ and $q_0q_k$,
 respectively. If $C_j$ and $C_k$ are disjoint, then, by Definition \ref{d:gates}  a
  gate  bigger than $R$ can be obtained  using these two cliques and the
  path described above depending on  cases $(a)$, $(b)$  or
  $(c)$. If $C_j$ and $C_k$ are non disjoint, then, by Lemma \ref{l:nodisj}  a
  gate  bigger than $R$ can be obtained  also using these two cliques and the
  path described above depending on  cases $(a)$, $(b)$ or
  $(c)$.  It contradicts the fact that $R$ is the biggest gate.
 \hfill$\Box$\end{pf}
Finally, to end the proof of Theorem \ref{t:bueno}, we will
describe below how to  obtain a new Helly $(d,2,2)$-representation
$\langle \mathcal{P'},T' \rangle$ of $G$ using a host tree $T'$
with fewer vertices of degree $d$. This contradicts the fact that
$\langle \mathcal{P},T \rangle$ is a representation minimizing the
number of vertices with degree $d$ and the proof follows.

To obtain $T'$ do:

subdivide the edge $q_0q_k$ of $T$ adding a new vertex
$\widetilde{q}_k$ adjacent to $q_0$ and to $q_k$;

for every  $k+1 \leq i \leq s$,    remove the edge $q_0q_{i}$ and
add the edge  $\widetilde{q}_kq_i$.

To obtain $\mathcal{P'}$ do:

If $P\in \mathcal{P}$ and there exist $k+1\leq i<j \leq s$ such
that $\{q_i,q_0,q_j\} \in  V(P)$, then replace $P$ by the path
$P'$ obtained from $P$ by removing the vertex $q_0$ and the edges
$q_0q_i$ and $q_0q_j$ and adding the vertex $\widetilde{q}_k$ and
the edges $\widetilde{q}_kq_i$ and $\widetilde{q}_kq_j$.

If $P\in \mathcal{P}$ and there exists $k+1\leq i \leq s$ such
that $\{q_i,q_0\} \in  V(P)$ and $P$ is not in the previous case
(then, by Claim \ref{c:q0qk} and the fact that $G_z$ is a
connected component of $G-A$, we have that $\{q_0,q_k\}\subseteq
V(P)$),  then replace $P$ by the path $P'$ obtained from $P$ by
removing the vertex $q_0$ and the edges $q_0q_i$ and $q_0q_k$ and
adding the vertex $\widetilde{q}_k$ and the edges
$\widetilde{q}_kq_i$ and $\widetilde{q}_kq_i$.

No other path is modified except for the fact of subdividing the
edge $q_0q_k$ if necessary.
 \hfill$\Box$\end{pf}
%
\begin{cor} \label{c:igual}  Helly $EPT \cap [h,2,2]=$ Helly $[h,2,2]$
 for any  $h\geq
    3$.\end{cor}
\begin{pf} Clearly, Helly $[h,2,2]\subseteq $ Helly $EPT \cap [h,2,2]$.

Assume, in order to derive a contradiction, that $G\in$  Helly
$EPT \cap [h,2,2]$ and $G\not\in $ Helly $[h,2,2]$. By
\ref{t:bueno}, $G$ contains a $k$-gate as induced subgraph for
some $k>h$. Thus by Theorem \ref{t:multipie}, any $EPT$
representation of $G$ contains a multipie of size $k$. This
contradicts the fact that $G\in [h,2,2]$. \hfill $\Box$\end{pf}
%
\section{Decomposition by clique separators and Complexity} \label{s:complex}
%
A clique $C$ of a connected graph $G$ is a \textbf{separator} if
$G-C$ (the subgraph induced by $V(G)\setminus C$) is not
connected. An \textbf{atom} is a connected graph with no
separators.  In \cite{Monma}, a graph is progressively decompose
by clique separators to obtain a \textbf{clique decomposition
tree} with each leaf node being associated with an atom of $G$ and
each internal node being associated with a clique separator of
$G$. The atoms of $G$ are invariants.  The clique decomposition
can be computed in polynomial time. Both $EPT$ graph and Helly
$EPT$ graphs are characterize by their clique decomposition tree.
The characterization leads to an efficient algorithm to recognize
Helly $EPT$ graphs but does not to recognize $EPT$ graphs.
%
%
\begin{lem} \label{l:atom-gate} If $H$ is a Helly $EPT$ atom with exactly $k\geq 4$ cliques
then $H$ has a $k$-gate as induced subgraph.
\end{lem}
%
\begin{pf}  Assume, in order to derive a contradiction, that $H$ has no
$k$-gates, then it has no $t$-gates for any $t\geq k$. Thus, by
Theorem \ref{t:bueno}, there exists $h\leq k-1$  such that $H\in$
Helly $[h,2,2]-$ Helly $[h-1,2,2]$. Let $\langle \mathcal{P},
T\rangle $ be a Helly $(h,2,2)$-representation of $H$ minimizing
the number of edges of the host tree $T$, this implies that $K_e$
is a clique of $H$ for every $e\in E(T)$, moreover $\mid E(T)
\mid=k$. On the other hand, since $H$ is an atom, $T$ must be a
star (otherwise there exists an edge $e$ of the host tree such
that $K_e$ is a cut clique). It follows that $h=k$, in
contradiction with the fact that $h< k$. \hfill $\Box$\end{pf}
%
\begin{lem} \label{l:gate-atom} Let $H$ be an $k$-gate. If $H$ is
an induced subgraph of a graph $G$, then $H$ is an induced
subgraph of some atom of $G$.\end{lem}
%
\begin{pf} It is enough to prove that a gate has no clique
separators which follows trivially from the recursive definition
of gates.
 \hfill $\Box$\end{pf}
\begin{thm} \label{t:HEPT} Let $G$ be a Helly $EPT$ graph and $h\geq 3$. Then,
$G\in$ Helly $[h,2,2]$ if an only if
every atom of $G$ has at most $h$ cliques.\end{thm}
%
\begin{pf} If $G\in $ Helly $[h,2,2]$ then, by Theorem
\ref{t:bueno}, $G$ has no gates of size grater than $h$ as induced
subgraphs. Thus, by Lemma \ref{l:atom-gate}, $G$ has no atoms with
more than $h$ cliques.

Conversely, assume, in order to obtain a contradiction, that
$G\not\in$ Helly $[h,2,2]$. Thus, by Theorem \ref{t:bueno}, $G$
has a $k$-gate $H$ as induced subgraph, for some $k>h$. By Lemma
\ref{l:gate-atom}, $H$ is an induced subgraph of some atom of $G$.
It implies that the atom has at least $k$ cliques, which
contradicts the assumption.\hfill $\Box$\end{pf}

We will consider the following two problems, the first is posed
for a given fixed  $h\geq 4$.

\begin{quote}
\textbf{RECOGNIZING HELLY} $\mathbf{[h,2,2]}$ \textbf{GRAPHS}\\
\underline{Input}: A connected graph $G$.\\
\underline{Question}: Does $G$ belong to Helly
$[h,2,2]$?\end{quote}

\begin{quote}
\textbf{CHEAPEST REPRESENTATION}\\
\underline{Input}: A connected  graph $G$.\\
\underline{Goal}: Determine the minimum $h\geq 2$ such that $G\in
$ Helly $[h,2,2]$. \end{quote}

Clearly an efficient solution of the latter implies an efficient
solution of the former.
%
\begin{thm} The problem CHEAPEST REPRESENTATION is polynomail times solvable.\end{thm}
    \begin{pf} Using the efficient algorithm described in  \cite{Monma}, determine
    whether
    the given graph $G$ belongs to Helly $EPT$ or not.
If it does then determine for each atom $G_i$ of $G$ its  number
of cliques, say $k_i$. Notice that it can be done efficiently
since
 the total number of cliques of a Helly
$EPT$ graphs $G$ is at most $\lfloor\frac{ 3 \mid V(G)\mid -
4}{2}\rfloor$  \cite{Monma}. Let $k$ be the maximum $k_i$.

 If
$k\leq 3$, then every atom is chordal which implies $G\in Chordal
\cap EPT=[3,2,2]$ (see \cite{Monma} and \cite{B4}). Now test
whether  $G$ is an interval graph or not and answer $h=2$ in an
affirmative case and $h=3$ otherwise.

If $k\geq 4$, by Theorem \ref{t:HEPT}, $G\in$ Helly $[k,2,2]$ and
$G\not\in$ Helly $[k-1,2,2]$, thus let  $h=k$.
 \hfill$\Box$\end{pf}

\end{document}